\documentclass[reqno,a4paper]{amsart}%
\usepackage{url,amssymb}
\usepackage[hypertex]{hyperref}
\allowdisplaybreaks[4]
\theoremstyle{plain}
\newtheorem{thm}{Theorem}
\newtheorem{prop}{Proposition}
\theoremstyle{remark}
\newtheorem{rem}{Remark}
\DeclareMathOperator{\td}{d\mspace{-1mu}}
\DeclareMathOperator{\re}{Re}

\newcommand{\tn}{\mathbb{N}}
\newcommand{\tr}{\mathbb{R}}
\date{October 29, 2003; Revised on April 22, 2004 at The Hong Kong University; Revised again in July, 2005}
\date{}

\begin{document}

\title[Some logarithmically completely monotonic functions]
{Some logarithmically completely monotonic functions related to the gamma function}

\author[F. Qi]{Feng Qi}
\address[F. Qi]{Research Institute of Mathematical Inequality Theory, Henan Polytechnic University, Jiaozuo City, Henan Province, 454010, China}
\email{\href{mailto: F. Qi <qifeng618@gmail.com>}{qifeng618@gmail.com}, \href{mailto: F. Qi <qifeng618@hotmail.com>}{qifeng618@hotmail.com}, \href{mailto: F. Qi <qifeng618@qq.com>}{qifeng618@qq.com}}
\urladdr{\url{http://qifeng618.spaces.live.com}}

\author[B.-N. Guo]{Bai-Ni Guo}
\address[B.-N. Guo]{School of Mathematics and Informatics,
Henan Polytechnic University, Jiaozuo City, Henan Province, 454010,
China} \email{\href{mailto: B.-N. Guo
<bai.ni.guo@gmail.com>}{bai.ni.guo@gmail.com}, \href{mailto: B.-N.
Guo <bai.ni.guo@hotmail.com>}{bai.ni.guo@hotmail.com}}
\urladdr{\url{http://guobaini.spaces.live.com}}

\begin{abstract}
In this article, logarithmically complete monotonicity properties of some
functions such as $\frac1{[\Gamma(x+1)]^{1/x}}$, $\frac{[{\Gamma(x+\alpha+1)}]^{1/(x+\alpha)}}{[{\Gamma(x+1)}]^{1/x}}$, $\frac{[\Gamma(x+1)]^{1/x}}{(x+1)^\alpha}$ and $\frac{[\Gamma(x+1)]^{1/x}}{x^\alpha}$ defined in $(-1,\infty)$ or $(0,\infty)$ for given real number $\alpha\in\tr$ are obtained, some known results are recovered, extended and generalized. Moreover, some basic properties of the logarithmically completely monotonic functions are established.
\end{abstract}



\thanks{This paper was typeset using \AmS-\LaTeX}

\maketitle

\section{Introduction}\label{minussec1}

Recall \cite[Chapter~XIII]{mpf} and \cite[Chapter~IV]{widder} that a function $f$ is said to be completely monotonic on an interval $I$ if $f$ has derivatives of all orders on $I$ and
\begin{equation}
(-1)^{k}f^{(k)}(x)\ge0
\end{equation}
for all $k\geq0$ on $I$. For our own convenience, let $\mathcal{C}[I]$ denote the set of completely monotonic functions on $I$. The well-known Bernstein's Theorem in \cite[p.~160, Theorem~12a]{widder} states that a function $f$ on $[0,\infty)$ is completely monotonic if and only if there exists a bounded and non-decreasing function $\alpha(t)$ such that
\begin{equation}
f(x)= \int_0^\infty e^{-xt}\td\alpha(t)
\end{equation}
converges for $x\in[0,\infty)$. This tells us that
$f\in\mathcal{C}[[0,\infty)]$ if and only if it is a Laplace
transform of the measure $\alpha$. There have been a lot of
literature about the completely monotonic functions, for examples,
\cite{Alzer, alzer-berg, alzbergii, grin-ismail, a27, mpf,
pams-62-rgmia, pams-62, cubo, cubo-rgmia, compmon, qx2, haerc1,
haer, haerc2, vogt, widder} and references therein.
\par
Recall also \cite{Atanassov, compmon2, minus-one-rgmia} that a positive function
$f$ is said to be logarithmically completely monotonic on an interval $I$ if $f$
has derivatives  of all orders on $I$ and
\begin{equation}\label{lcmf-ineq}
(-1)^n[\ln f(x)]^{(n)}\ge 0
\end{equation}
for all $ x\in I$ and $n\in\mathbb{N}$. For simplicity, let $\mathcal{L}[I]$ stand for the set of logarithmically completely monotonic functions on $I$.
\par
Among other things, it is proved in \cite{CBerg, bochner, compmon2, minus-one-rgmia, Mon-Two-Seq-complete.tex} that a logarithmically completely
monotonic function is always completely monotonic, that is,
$\mathcal{L}[I]\subset\mathcal{C}[I]$, but not conversely, since a convex
function may not be logarithmically convex (see \cite[p.~7,
Remark.~1.16]{cposa}).
\par
Recall \cite{widder} that a function $f$ defined in $(0,\infty)$ is called a Stieltjes transform if it can be of the form
\begin{equation}
f(x)=a+\int_0^\infty\frac1{s+x}{\td\mu(s)},
\end{equation}
where $a$ is a nonnegative number and $\mu$ a nonnegative measure on $[0,\infty)$ satisfying
\begin{equation}
\int_0^\infty\frac1{1+s}\td\mu(s)<\infty.
\end{equation}
The set of Stieltjes transforms is denoted by $\mathcal{S}$.
\par
Motivated by the papers \cite{minus-one-rgmia,auscmrgmia}, among other things, it
is further revealed in \cite{CBerg} that
\begin{equation}
\mathcal{S}\setminus\{0\}\subset\mathcal{L}[(0,\infty)]\subset\mathcal{C}[(0,\infty)].
\end{equation}
In \cite[Theorem~1.1]{CBerg} and \cite{grin-ismail,auscm} it is pointed out
that logarithmically completely monotonic functions on $(0,\infty)$ can be
characterized as the infinitely divisible completely monotonic functions
studied by Horn in \cite[Theorem~4.4]{horn}. The functions in $\mathcal{L}[I]$
are also characterized by $-\frac{f'}f\in\mathcal{C}[I]$. Recently it is found
that a finer inclusion
\begin{equation}
\mathcal{S}\subset\mathcal{C}^\ast[(0,\infty)] \subset\mathcal{L}[(0,\infty)]\subset\mathcal{C}[(0,\infty)]
\end{equation}
had been established in \cite[Section~14.2, pp.~122--127]{berg-forst} and \cite{haer}, where $\mathcal{C}^\ast[(0,\infty)]$ denotes the set
\begin{equation}
\biggl\{f(x)\,\bigg|\,\biggl[\frac1{f(x)}\biggr]'\in \mathcal{C}[(0,\infty)]\biggr\}.
\end{equation}
In \cite[p.~122]{berg-forst}, it was proved that \begin{equation}
\mathcal{H}\setminus\{0\}\subset\mathcal{P}=\mathcal{C}^\ast[(0,\infty)]
\end{equation}
and it was told that this is a theorem of F. Hirsch with due reference. This result says that if $f>0$ and $f\in\mathcal{H}$ then $\frac1f$ is the Laplace transform of a potential kernel, hence the Laplace transform of an infinitely divisible measure and $f$ is a Bernstein function, i.e. a positive function whose derivative is completely monotonic. On \cite[p.~127]{berg-forst} it is proved that $\mathcal{S}\subset\mathcal{H}$.
\par
From Bernstein's Theorem it also follows that completely monotonic functions on
$(0,\infty)$ are always strictly completely monotonic unless they are constant,
see \cite{dub,auscm} and \cite[p.~11]{haer}. Also it follows that a
logarithmically completely monotonic function on $(0,\infty)$ is strictly so
unless it is of the form $c\exp(-\alpha x)$ for $c>0$ and $\alpha\ge0$, so there
is no need to discuss the sharpening with ``strictly" in general. If its representing measure of a function
$f$ as a Stieltjes transform is concentrated on $[a,\infty)$ with $a>0$, then
$f\in\mathcal{L}[(-a,\infty)]$.
\par
The classical Euler gamma function is usually defined  for $\re z>0$ by
\begin{equation}
\Gamma(z)=\int^\infty_0t^{z-1} e^{-t}\td t.
\end{equation}
The logarithmic derivative of the gamma function
\begin{equation}
\psi(x)=\frac{\Gamma'(x)}{\Gamma(x)}
\end{equation}
is called the psi or digamma function and $\psi^{(n)}(x)$ for
$n\in\tn$ the polygamma functions. It is well-known that the gamma
function is a very important classical special function and has many
applications \cite{abram,bochner,emot,magnus}. One of the reasons
why the gamma function is still interesting, although nearly three
centuries have elapsed after its first appearance, is that it has
many applications to various areas of mathematics ranging from
probability theory to number theory and function theory.
(Logarithmically) completely monotonic functions have applications
in many branches. For example, they play a role in complex analysis,
number theory, potential theory, probability theory \cite{bochner},
physics \cite{magnus}, numerical and asymptotic analysis, integral
transforms \cite{widder}, and combinatorics. Some related references
are listed in
\cite{Alzer,alzer-berg,alzbergii,CBerg,grin-ismail,mpf,rsg,indon,auscm,widder}.
In recent years, inequalities and (logarithmically) completely
monotonic functions involving the gamma, psi, or polygamma functions
are established by some mathematicians (see
\cite{Allasia-Gior-Pecaric-MIA-02, Alzer, alzer-berg, alzbergii,
tamk05364, chenqig, jmaa05726, glob051, chenqigam, chenqigam-rgmia,
clark-ismail, grin-ismail, qigtai, kershawl, miller, mrigf-ii,
mrigf-ii-rgmia, qccg, geom1, notes-best, notes-best-rgmia, vogt} and
the references therein).
\par
In this paper, using Leibniz's Identity, the discrete and integral
representations of polygamma functions and other analytic techniques, some
functions such as
\begin{align}
&\frac1{[\Gamma(x+1)]^{1/x}},&&
\frac{[{\Gamma(x+\alpha+1)}]^{1/(x+\alpha)}}{[{\Gamma(x+1)}]^{1/x}},&&
\frac{[\Gamma(x+1)]^{1/x}}{(x+1)^\alpha}&& \text{and}&&
\frac{[\Gamma(x+1)]^{1/x}}{x^\alpha}
\end{align}
with $x\in(-1,\infty)$ or $x\in(0,\infty)$ for given real number $\alpha\in\tr$ are shown to be logarithmically completely monotonic. Moreover, some basic properties of the logarithmically completely monotonic functions are established.
\par
Our main results are as follows.

\begin{thm}\label{gammc}
\begin{equation}
\frac1{[\Gamma(x+1)]^{1/x}}\in\mathcal{L}[(-1,\infty)].
\end{equation}
\end{thm}

\begin{thm}\label{gammd}
\begin{equation}
\frac{[\Gamma(x+1)]^{1/x}}{(x+1)^\alpha}\in\mathcal{L}[(-1,\infty)]
\end{equation}
if and only if $\alpha\ge1$.
\end{thm}

\begin{thm}\label{gamme}
\begin{equation}
\frac{x^\alpha}{[\Gamma(x+1)]^{1/x}}\in\mathcal{L}[(0,\infty)]
\end{equation}
if and only if $\alpha\le0$.
\begin{equation}
\frac{[\Gamma(x+1)]^{1/x}}{x^\alpha}\in\mathcal{L}[(0,\infty)]
\end{equation}
if and only if $\alpha\ge1$.
\end{thm}

\begin{thm}\label{thmadd}
Let
\begin{equation}\label{taudfn}
\tau(s,t)=\frac1s\biggl[t-(t+s+1)\Bigl(\frac{t}{t+1}\Bigr)^{s+1}\biggr]
\end{equation}
for $(s,t)\in(0,\infty)\times(0,\infty)$ and $\tau_0=\tau(s_0,t_0)>0$ be the
maximum of $\tau(s,t)$ on the set $\tn\times(0,\infty)$. Then for any given
real number $\alpha$ satisfying $\alpha\le\frac1{1+\tau_0}<1$,
\begin{equation}
\frac{(x+1)^\alpha}{[\Gamma(x+1)]^{1/x}}\in\mathcal{L}[(-1,\infty)].
\end{equation}
\end{thm}

\begin{thm}\label{thm:NAME5}
For $\alpha\le0$ such that $x^\alpha$ is real in $(-1,0)$,
\begin{equation}
\frac{[\Gamma(x+1)]^{1/x}}{x^\alpha}\in\mathcal{L}[(-1,0)].
\end{equation}
For $\alpha\ge1$ such that $x^\alpha$ is real in $(-1,0)$,
\begin{equation}
\frac{x^\alpha}{[\Gamma(x+1)]^{1/x}}\in\mathcal{L}[(-1,0)].
\end{equation}
\end{thm}

As basic properties of the logarithmically completely monotonic functions, we obtain
the following theorems.

\begin{thm}\label{ratthm}
Let $f(x)\in\mathcal{L}[I]$. Then $\frac{f(x)}{f(x+\alpha)}\in\mathcal{L}[J]$
if and only if $\alpha>0$, where $J=I\cap\{x+\alpha\in I\}$.
\end{thm}

\begin{thm}
\label{thm7minus} Let $f_i(x)\in\mathcal{L}[I]$ and $\alpha_i\ge0$ for $1\le i\le n$
with $n\in\tn$. Then
\begin{equation}
\prod_{i=1}^n[f_i(x)]^{\alpha_i}\in\mathcal{L}[I].
\end{equation}
\end{thm}

\begin{thm}\label{compositelogconv}
Let $h'(x)\in\mathcal{C}[I]$ and $f(x)\in\mathcal{L}[h(I)]$. Then $f\circ
h(x)=f\big(h(x)\big)\in\mathcal{L}[I]$.
\end{thm}

In Section \ref{sec:Proofs-theorem}, we are about to give proofs of these theorems.
In Section \ref{cores}, some remarks are given, some new results are deduced, and
some known results are recovered, as applications of these theorems.

\section{Proofs of theorems}\label{sec:Proofs-theorem}

It is well-known (see \cite{abram,emot,wange,wang} and \cite[p.~16]{magnus}) that the
polygamma functions $\psi^{(k)}(x)$ can be expressed for $x>0$ and $k\in\mathbb{N}$
as
\begin{align}
\label{difpsi} \psi^{(k)}(x)&=(-1)^{k+1}k!
\sum_{i=0}^{\infty}\frac{1}{(x+i)^{k+1}}\intertext{or}
\psi
^{(k)}(x)&=(-1)^{k+1}\int_{0}^{\infty}\frac{t^{k}e^{-xt}}{1-e^{-t}}\td t.
\label{psim}
\end{align}

\begin{proof}[The first proof of Theorem \ref{gammc}]
Let
\begin{equation}\label{gxdef}
g(x)=\begin{cases}
\dfrac{\ln\Gamma(x+1)}{x},&x\ne0\\
-\gamma,&x=0
\end{cases}
\end{equation}
for $x\in(-1,\infty)$, where $\gamma=0.57721566\dotsm$ is the Euler-Mascheroni
constant. By direct calculation and using Leibniz's Identity, we obtain for
$n\in\tn$,
\begin{align}
g^{(n)}(x)&=\begin{cases}
\dfrac{1}{x^{n+1}}\displaystyle\sum_{k=0}^n\dfrac{(-1)^{n-k}n!x^k\psi^{(k-1)}(x+1)}{k!}
\triangleq \dfrac{h_n(x)}{x^{n+1}},&x\ne0,\\
\dfrac{\psi^{(n)}(1)}{n+1},&x=0,
\end{cases}\\
h_n'(x)&=x^n\psi^{(n)}(x+1)\begin{cases}
>0&\text{in $(0,\infty)$ if $n$ is odd},\\
\le0&\text{in $(-1,0]$ if $n$ is odd},\\
\le0&\text{in $(-1,\infty)$ if $n$ is even},
\end{cases}
\end{align}
where
$$
\psi^{(-1)}(x+1)=\ln\Gamma(x+1)\quad \text{and} \quad\psi^{(0)}(x+1)=\psi(x+1).
$$
Hence, if $n$ is odd the function $h_n(x)$ increases in $(0,\infty)$ and
decreases in $(-1,0)$, if $n$ is even it decreases in $(-1,\infty)$. Since
$h_n(0)=0$, it is easy to see that $h_n(x)\ge0$ in $(-1,\infty)$ if $n$ is odd
and that $h_n(x)\ge0$ in $(-1,0)$ and $h_n(x)\le0$ in $(0,\infty)$ if $n$ is
even. Thus, in the interval $(-1,\infty)$, the function $g^{(n)}(x)\ge0$ if
$n$ is odd and $g^{(n)}(x)\le0$ if $n$ is even. Since
$$
\lim_{x\to\infty}\frac{\psi^{(k)}(x+1)}{x^{n-k}}=0
$$
for $-1\le k\le n-1$, it
follows that $\lim_{x\to\infty}g^{(n)}(x)=0$. Consequently,
$$
(-1)^{n+1}g^{(n)}(x)>0
$$
in $(-1,\infty)$ for $n\in\tn$. This implies that
$$
(-1)^{k}\{\ln [\Gamma(x+1)]^{1/x}\}^{(k)}<0
$$
in $(-1,\infty)$ for $k\in\tn$ and the function $\frac{1}{[\Gamma(x+1)]^{1/x}}$ is logarithmically completely
monotonic in $(-1,\infty)$.
\end{proof}

\begin{proof}[The second proof of Theorem \ref{gammc}]
It is not difficult to see that
\begin{gather}
g(x)=\frac{\ln\Gamma(x+1)-\ln\Gamma(1)}{x} =\frac1x\int_0^x\psi(t+1)\td t
=\int_0^1\psi(xs+1)\td s \intertext{and}
g^{(n)}(x)=\int_0^1s^n\psi^{(n)}(xs+1)\td s.\label{gnint}
\end{gather}
Thus, the required result follows from using formula \eqref{difpsi} or \eqref{psim}
in \eqref{gnint}.
\end{proof}

\begin{proof}[Proof of Theorem \ref{gammd}]
Let
\begin{equation}\label{nudfn}
\nu_\alpha(x)=\begin{cases}
\dfrac{[\Gamma(x+1)]^{1/x}}{(x+1)^\alpha},&x\ne0\\e^{-\gamma},&x=0
\end{cases}
\end{equation}
for $x\in(-1,\infty)$. Then for $n\in\tn$, by using \eqref{difpsi},
\begin{align}
\ln \nu_\alpha(x)&=\frac{\ln\Gamma(x+1)}x-\alpha\ln(x+1),\\
[\ln \nu_\alpha(x)]^{(n)}&=\frac1{x^{n+1}}\biggl[h_n(x)+\frac{(-1)^{n}(n-1)!\alpha
x^{n+1}}{(x+1)^n}\biggr]\triangleq\frac{\mu_{\alpha,n}(x)}{x^{n+1}},
\end{align}
\begin{align}\notag
\mu_{\alpha,n}'(x)&=x^n\psi^{(n)}(x+1)+\frac{(-1)^{n}(n-1)!\alpha
x^{n}(x+n+1)}{(x+1)^{n+1}}\\\notag
&=x^n\biggl[\psi^{(n)}(x+1)+\frac{(-1)^{n}(n-1)!\alpha}
{(x+1)^{n}}+\frac{(-1)^{n}n!\alpha}{(x+1)^{n+1}}\biggr]\\\notag
&=x^n\biggl\{(-1)^{n+1}n!\sum_{i=1}^{\infty}\frac{1}{(x+i)^{n+1}}\\\notag
&\quad+{(-1)^{n}(n-1)!\alpha}
\sum_{i=1}^\infty\biggl[\frac1{(x+i)^{n}}-\frac1{(x+i+1)^{n}}\biggr]\label{minus-one-32}\\
&\quad+{(-1)^{n}n!\alpha}\sum_{i=1}^\infty
\biggl[\frac1{(x+i)^{n+1}}-\frac1{(x+i+1)^{n+1}}\biggr]\biggr\}\\\notag
&=(-1)^n(n-1)!x^n\sum_{i=1}^\infty
\biggl[\frac{\alpha}{(x+i)^n}-\frac{\alpha}{(x+i+1)^n}\\*\notag
&\quad-\frac{n\alpha}{(x+i+1)^{n+1}}+\frac{n(\alpha-1)}{(x+i)^{n+1}}\biggr]\notag
\end{align}
\begin{align}\notag
&=(n-1)!(-x)^n\sum_{i=1}^\infty\frac{[\alpha
y+n(\alpha-1)](y+1)^{n+1}-\alpha(y+n+1)y^{n+1}}{y^{n+1}(y+1)^{n+1}}\\\notag
&=(n-1)!(-x)^n\sum_{i=1}^\infty
\frac{\alpha[(y+n)(y+1)^{n+1}-(y+n+1)y^{n+1}]-n(y+1)^{n+1}}{y^{n+1}(y+1)^{n+1}}\\\notag
&=n!(-x)^n\sum_{i=1}^\infty\frac{1}{y^{n+1}}
\biggl\{\alpha\biggl[1+\frac1n\biggl\langle
y-(y+n+1)\biggl(\frac{y}{y+1}\biggr)^{n+1}\biggr\rangle\biggr]-1\biggr\},
\end{align}
where $y=x+i>0$.
\par
In \cite[p.~28]{bullen}, \cite[p.~154]{kuangb} and \cite{3rded}, Bernoulli's
inequality states that if $x\ge-1$ and $x\ne0$ and if $\alpha>1$ or if
$\alpha<0$ then
$$
(1+x)^\alpha>1+\alpha x.
$$
This means that
$$
1+\frac{s+1}t<\biggl(1+\frac1t\biggr)^{s+1}
$$
which is equivalent to
$$
t-(t+s+1)\biggl(\frac{t}{t+1}\biggr)^{s+1}>0
$$
for $s>0$ and $t>0$, then the
function $\tau(s,t)$ defined by \eqref{taudfn} is positive for
$(s,t)\in(0,\infty)\times(0,\infty)$.
\par
Since $\tau(s,t)>0$, it is deduced that $$[\alpha
y+n(\alpha-1)](y+1)^{n+1}-\alpha(y+n+1)y^{n+1}>0$$ for $y=x+i>0$ and $n\in\tn$ if
$\alpha\ge1$. This means that for $\alpha\ge1$,
\begin{equation*}
\mu_{\alpha,n}'(x)
\begin{cases}
>0&\text{in $(-1,0)\cup(0,\infty)$ if $n$ is even},\\
>0&\text{in $(-1,0)$ if $n$ is odd},\\
<0&\text{in $(0,\infty)$ if $n$ is odd},
\end{cases}
\end{equation*}
hence, it is obtained that the function $\mu_{\alpha,n}(x)$ is strictly
increasing in $(-1,\infty)$ if $n$ is even and that the function
$\mu_{\alpha,n}(x)$ is strictly increasing in $(-1,0)$ and strictly decreasing
in $(0,\infty)$ if $n$ is odd. Since $\mu_{\alpha,n}(0)=0$, it follows that
$\mu_{\alpha,n}(x)\le0$ in $(-1,\infty)$ if $n$ is odd and that
$\mu_{\alpha,n}(x)\le0$ in $(-1,0)$ and $\mu_{\alpha,n}(x)>0$ in $(0,\infty)$
if $n$ is even. From $\lim_{x\to\infty}[\ln \nu_\alpha(x)]^{(n)}=0$, it is
concluded that $[\ln \nu_\alpha(x)]^{(n)}\ge0$ in $(-1,\infty)$ if $n$ is even
and $[\ln \nu_\alpha(x)]^{(n)}\le0$ in $(-1,\infty)$ if $n$ is odd, which is
equivalent to $(-1)^n[\ln \nu_\alpha(x)]^{(n)}>0$ in $x\in(-1,\infty)$ for
$n\in\tn$ and $\alpha\ge1$. Hence, if $\alpha\ge1$, the function
$\frac{[\Gamma(x+1)]^{1/x}}{(x+1)^\alpha}$ is logarithmically completely
monotonic in $(-1,\infty)$.
\par
Conversely, if the function $\frac{[\Gamma(x+1)]^{1/x}}{(x+1)^\alpha}$ is
logarithmically completely monotonic in $(-1,\infty)$, then $[\ln
\nu_\alpha(x)]'\le0$ which is equivalent to
\begin{align*}
\alpha&\ge\frac{x+1}{x^2}[x\psi(x+1)-\ln\Gamma(x+1)]\\
&=\bigg(1+\frac1x\bigg)\bigg[\psi(x+1)-\frac{\ln\Gamma(x+1)}x\bigg]\\
&\to\bigg(1+\frac1x\bigg)\bigg\{\ln(x+1)-\frac1{2(x+1)}-\frac1{12(x+1)^2}+O\bigg(\frac1{x+1}\bigg)\\
&\quad-\frac1x\bigg[\bigg(x+\frac12\bigg)\ln(x+1)-x-1+\frac{\ln(2\pi)}2
+\frac1{12(x+1)}+O\bigg(\frac1{x+1}\bigg)\bigg]\bigg\}\\* &\to1
\end{align*}
as $x\to\infty$ by using the following formulas (see \cite{abram,magnus,wange,wang})
\begin{align}
\ln\Gamma(x)&=\bigg(x-\frac12\bigg)\ln x-x+\frac{\ln(2\pi)}{2}
+\frac1{12x}+O\bigg(\frac1x\bigg)\label{psisymp1}\intertext{and}
\psi(x)&=\ln
x-\frac1{2x}-\frac1{12x^2}+O\bigg(\frac1{x^2}\bigg)\label{psisymp}
\end{align}
as $x\to\infty$.
\end{proof}

\begin{proof}[Proof of Theorem \ref{gamme}]
If $\alpha\le0$, the logarithmically complete monotonicity of the function
$\frac{x^\alpha}{[\Gamma(x+1)]^{1/x}}$ in $(0,\infty)$ follows from the
similar arguments as in the proofs of Theorem~\ref{gammd}.
\par
If $\frac{x^\alpha}{[\Gamma(x+1)]^{1/x}}$ in $(0,\infty)$ is logarithmically
completely monotonic, then its logarithmic derivative
$$
\frac{\alpha}x+\frac{\ln\Gamma(1+x)-{x\psi(1+x)}}{x^2}
$$
is negative in $(0,\infty)$. Since
\begin{equation}
\label{nbx}
\lim_{x\to0^+}\frac{\ln\Gamma(1+x)-{x\psi(1+x)}}{x^2}=-\frac{\pi^2}{12}
\end{equation}
by L'Hospital rule and
\begin{equation}
\label{asb} \lim_{x\to0^+}\frac{\alpha}{x}=\begin{cases}
\infty,& \text{if $\alpha>0$},\\
0,&\text{if $\alpha=0$},\\
-\infty,&\text{if $\alpha<0$},
\end{cases}
\end{equation}
then it must hold that $\alpha\le0$.
\par
The rest proofs of Theorem \ref{gamme} are similar to the proofs of
Theorem~\ref{gammd}, so we omit it.
\end{proof}

\begin{proof}[Proof of Theorem \ref{thmadd}]
Since $\tau(s,t)>0$, which has been proved in Theorem \ref{gammd} by utilizing
Bernoulli's inequality, it is clear that $\tau_0>0$. When
$\alpha\le\frac1{1+\tau_0}<1$, from \eqref{minus-one-32} it follows that
$\mu_{\alpha,n}'(x)\le0$ and $\mu_{\alpha,n}(x)$ is decreasing in
$(-1,\infty)$ if $n$ an even integer and that $\mu_{\alpha,n}'(x)\le0$ and
$\mu_{\alpha,n}(x)$ is decreasing in $(-1,0)$ and $\mu_{\alpha,n}'(x)\ge0$ and
$\mu_{\alpha,n}(x)$ is increasing in $(0,\infty)$ if $n$ an odd integer. Since
$\mu_{\alpha,n}(0)=0$ and $\lim_{x\to\infty}[\ln \nu_\alpha(x)]^{(n)}=0$, we
have $[\ln\nu_\alpha(x)]^{(n)}<0$ for $n$ being an even and
$[\ln\nu_\alpha(x)]^{(n)}>0$ for $n$ being an odd in $(-1,\infty)$, this
implies that $(-1)^{n+1}[\ln\nu_\alpha(x)]^{(n)}>0$ in $(-1,\infty)$ for
$n\in\tn$. Therefore $\nu_\alpha(x)$ is strictly increasing and
$(-1)^{n-1}\{[\ln\nu_\alpha(x)]'\}^{(n-1)}>0$ in $(-1,\infty)$ for $n\in\tn$.
Hence, if $\alpha\le\frac1{1+\tau_0}$, then the function
$\frac{(x+1)^\alpha}{[\Gamma(x+1)]^{1/x}}$ is logarithmically completely
monotonic in $(-1,\infty)$. The proof of Theorem \ref{thmadd} is complete.
\end{proof}

\begin{proof}[Proof of Theorem \ref{thm:NAME5}]
This follows from modified arguments of above theorems.
\end{proof}

\begin{proof}[Proof of Theorem \ref{ratthm}]
Let $\mathcal{F}_\alpha(x)=\frac{f(x)}{f(x+\alpha)}$ for $\alpha>0$. Since
$f(x)$ is logarithmically completely monotonic, by definition we have
$(-1)^k[\ln f(x)]^{(k)}\ge0$ for $k\in\tn$, which is equivalent to $[\ln
f(x)]^{(2i)}\ge0$ and $[\ln f(x)]^{(2i-1)}\le0$ for $i\in\tn$, and $[\ln
f(x)]^{(2i)}$ is decreasing and $[\ln f(x)]^{(2i-1)}$ is increasing. So
$$
[\ln\mathcal{F}_\alpha(x)]^{(2i)}=[\ln f(x)]^{(2i)}-[\ln f(x+\alpha)]^{(2i)}\ge0
$$
and $[\ln\mathcal{F}_\alpha(x)]^{(2i-1)}\le0$ for
$\alpha>0$ and $i\in\tn$. The proof of Theorem \ref{ratthm} is complete.
\end{proof}

\begin{proof}[Proof of Theorem \ref{thm7minus}]
Let
$$
F_n(x)=\prod_{i=1}^n[f_i(x)]^{\alpha_i}.
$$
Then
$$
\ln F_n(x)=\sum_{i=1}^n\alpha_i\ln f_i(x)
$$
and
\begin{equation*}
(-1)^k[\ln F_n(x)]^{(k)}=\sum_{i=1}^n\alpha_i(-1)^k[\ln f_i(x)]^{(k)}
\end{equation*}
for $k\in\tn$. Since $f_i(x)\in\mathcal{L}[I]$, that is, $(-1)^k[\ln
f_i(x)]^{(k)}\ge0$, and $\alpha_i\ge0$, it is easy to see that $(-1)^k[\ln
F_n(x)]^{(k)}\ge0$ for $k\in\tn$. The proof of Theorem \ref{thm7minus} is
complete.
\end{proof}

\begin{proof}[Proof of Theorem \ref{compositelogconv}]
In \cite[No.~0.430.1]{grads} the formula for the $n$-th derivative of a composite
function is given by
\begin{gather}
\frac{\td^n}{\td
x^n}\big[f\big(h(x)\big)\big]=\sum_{k=1}^n\frac1{k!}f^{(k)}\big(h(x)\big)U_k(x),
\intertext{where} U_k(x)=\sum_{i=0}^{k-1}(-1)^i\binom{k}{i}[h(x)]^i \frac{\td^n}{\td
x^n}[h(x)]^{k-i}.
\end{gather}
From this it is deduced that $-\frac{f'(h(x))}{f(h(x))}\in\mathcal{C}[I]$
since $-\frac{f'(x)}{f(x)}$ is a completely monotonic function, which is
equivalent to $f(x)$ being logarithmically completely monotonic, and
$h'(x)\in\mathcal{C}[I]$. Therefore,
$(-1)^i\big[\frac{f'(h(x))}{f(h(x))}\big]^{(i)}\le0$ on the interval $I$ for
nonnegative integer $i$.
\par
Since $h'(x)\in\mathcal{C}[I]$, it is obtained that $(-1)^ih^{(i+1)}(x)\ge0$ on the
interval $I$ for nonnegative integer $i$.
\par
Hence, for $k\in\tn$,
\begin{multline*}
(-1)^k\big[\ln f\big(h(x)\big)\big]^{(k)}=(-1)^k\bigg[\frac{f'\big(h(x)\big)}
{f\big(h(x)\big)}h'(x)\bigg]^{(k-1)}\\
=(-1)^k\sum_{i=0}^{k-1}\binom{k-1}{i}\bigg[\frac{f'\big(h(x)\big)}
{f\big(h(x)\big)}\bigg]^{(i)}h^{(k-i)}(x)\\
=\sum_{i=0}^{k-1}\binom{k-1}{i}\bigg\{(-1)^i\bigg[\frac{f'\big(h(x)\big)}
{f\big(h(x)\big)}\bigg]^{(i)}\bigg\}\big[(-1)^{k-i}h^{(k-i)}(x)\big]\ge0.
\end{multline*}
The proof of Theorem \ref{compositelogconv} is complete.
\end{proof}

\section{Remarks and applications of theorems}\label{cores}

\begin{rem}
As said in \cite{CBerg} and done in various papers, the complete monotonicity
for special functions has been established by proving the stronger statement
that the function is logarithmically completely monotonic or is a Stieltjes
transform. In some concrete cases it is often easier to establish that a
function is logarithmically completely monotonic or is a Stieltjes transform
than to verify directly the complete monotonicity. One of the important values
of this paper might be owning to the standard or elementary proofs of some
theorems in this paper.
\end{rem}

\begin{rem}
It is remarked that many complete monotonicity results in \cite{Alzer, alzer-berg, alzbergii, clark-ismail, miller, vogt} and the references therein can be restated in terms of the logarithmically complete monotonicity indeed.
\end{rem}

\begin{rem}\label{rem:NAME}
In \cite{alzer-berg} and \cite[p.~83]{bochner}, the following result was
given: Let $f$ and $g$ be functions such that $f\circ g$ is defined. If
$f\in\mathcal{C}[(0,\infty)]$ and $g'\in\mathcal{C}[(0,\infty)]$, then $f\circ
g\in\mathcal{C}[(0,\infty)]$. Since the exponential function
$e^{-x}\in\mathcal{C}[(-\infty,\infty)]$, hence
$\mathcal{L}[(0,\infty)]\subset\mathcal{C}[(0,\infty)]$, a logarithmically
completely monotonic function is also completely monotonic. This gives an
alternative proof of the conclusion
$\mathcal{L}[(0,\infty)]\subset\mathcal{C}[(0,\infty)]$.
\end{rem}

\begin{rem}
In \cite{alzbergii, CBerg} it is shown that $\frac{[\Gamma(1+x)]^{1/x}}x\in\mathcal{S}$ and $\frac1{[\Gamma(1+x)]^{1/x}}\in\mathcal{S}$. Although these results are stronger than Theorem \ref{gammc} and parts of Theorem~\ref{gamme}, but the ranges of $x$ are extended to $(-1,\infty)$ in Theorem \ref{gammc} and a parameter $\alpha$ is considered in Theorem~\ref{gamme}. Consequently, Theorem \ref{gammc} and Theorem~\ref{gamme} still make sense.
\end{rem}

\begin{rem}
It is noted that a non-elementary argument for the sufficient part of
Theorem~\ref{gammd} was provided by an anonymous referee of this paper as
follows. Looking at what is really written in \cite{alzbergii}, which builds
on a technique from \cite{berg-pedersen-rocky}, it is easy to obtain that
\begin{equation*}
h(z)=\frac{\ln\Gamma(z+1)}z-\alpha\ln(z+1)=c+\int_1^\infty\biggl(\frac1{t+z}
-\frac{t}{1+t^2}\biggr)[\alpha-M(t)]\td t
\end{equation*}
with
$$
c=-\gamma+\sum_{k=1}^\infty\biggl(\frac1k-\arctan\frac1k\biggr)
$$
and $M(t)=\frac{k-1}t$ for $t\in(k-1,k]$ and $k=2,3,\dotsc$. For $\alpha\ge1$ one has $\alpha\ge M(t)$. Accordingly,
\begin{equation}
h'(t)=-\int_1^\infty\frac{\alpha-M(t)}{(t+z)^2}\td t<0,\qquad z>-1,
\end{equation}
that is, the function $h$ is decreasing with $-h'\in\mathcal{C}[(-1,\infty)]$,
which is the sufficient part of Theorem~\ref{gammd}.
\end{rem}

\begin{rem}
Now we give a weaker proof of Theorem \ref{gammd} and Theorem \ref{thmadd} by another
approach.
\par
It is well-known (see \cite{abram,wange,wang}) that for $x>0$ and $r>0$
\begin{equation}\label{fracint}
\frac1{x^r}=\frac1{\Gamma(r)}\int_0^\infty t^{r-1}e^{-xt}\td t.
\end{equation}
\par
Substituting \eqref{fracint} and \eqref{psim} into the second line of formula
\eqref{minus-one-32} yields
\begin{align*}
\mu_{\alpha,n}'(x)&=x^n\biggl[\psi^{(n)}(x+1)+\frac{(-1)^{n}(n-1)!\alpha}
{(x+1)^{n}}+\frac{(-1)^{n}n!\alpha}{(x+1)^{n+1}}\biggr]\\
&=x^n\bigg\{(-1)^{n+1}\int_0^\infty\frac{t^ne^{-(x+1)t}}{1-e^{-1}}\td t\\
&\quad+(-1)^n\alpha\bigg[\int_0^\infty t^{n-1}e^{-(x+1)t}\td t+\int_0^\infty
t^{n}e^{-(x+1)t}\td t\bigg]\bigg\}\\
&=(-1)^nx^n\int_0^\infty\bigg[\alpha-\frac{et}{(e-1)(1+t)}\bigg]
(1+t)t^{n-1}e^{-(x+1)t}\td t.
\end{align*}
It is clear that the function $\frac{t}{1+t}$ is increasing in $[0,\infty)$
with $0\le \frac{t}{1+t}<1$. Thus, the function
$(-1)^n\frac{\mu_{\alpha,n}'(x)}{x^n}$ is non-positive for $\alpha\le0$ and
positive for $\alpha\ge\frac{e}{e-1}$ in $(-1,\infty)$. Then
\subsubsection*{\rm{(1)}}
For $n$ is even,
\begin{equation}
\mu_{\alpha,n}'(x)
\begin{cases}
\le0&\text{if $\alpha\le0$},\\>0&\text{if $\alpha\ge\frac{e}{e-1}$}.
\end{cases}
\end{equation}
and
\begin{equation}
\mu_{\alpha,n}(x)
\begin{cases}
\text{is decreasing if $\alpha\le0$},\\\text{is increasing if
$\alpha\ge\frac{e}{e-1}$},
\end{cases}
\end{equation}
hence, from $\mu_{\alpha,n}(0)=0$, it follows that
\begin{equation}
\begin{cases}\text{if $\alpha\le0$, }&\mu_{\alpha,n}(x)
\begin{cases}
>0,&x\in(-1,0),\\<0,&x\in(0,\infty),
\end{cases}\\[1em]
\text{if $\alpha\ge\dfrac{e}{e-1}$, }&\mu_{\alpha,n}(x)
\begin{cases}
<0,&x\in(-1,0),\\>0,&x\in(0,\infty),
\end{cases}
\end{cases}
\end{equation}
therefore,
\begin{equation}
[\ln\nu_\alpha(x)]^{(n)}=\frac{\mu_{\alpha,n}(x)}{x^{n+1}}
\begin{cases}
\le0&\text{if $\alpha\le0$,}\\\ge0&\text{if $\alpha\ge\frac{e}{e-1}$;}
\end{cases}
\end{equation}
\subsubsection*{\rm{(2)}}
For $n$ is odd,
\begin{equation*}
\frac{\mu_{\alpha,n}'(x)}{x^n}\begin{cases}\ge0&\text{if $\alpha\le0$},\\
<0&\text{if $\alpha\ge\frac{e}{e-1}$},
\end{cases}
\end{equation*}
consequently,
\begin{equation*}
\begin{cases}
\text{if
$\alpha\le0$,}&\mu_{\alpha,n}'(x)\begin{cases}\ge0&x\in(0,\infty),\\\le0&x\in(-1,0),
\end{cases}\\[1em]
\text{if
$\alpha\ge\dfrac{e}{e-1}$,}&\mu_{\alpha,n}'(x)\begin{cases}<0&x\in(0,\infty),\\>0&x\in(-1,0),
\end{cases}
\end{cases}
\end{equation*}
and
\begin{equation*}
\begin{cases}
\text{if $\alpha\le0$,}&\mu_{\alpha,n}(x)\begin{cases}\text{is increasing in
$(0,\infty)$,}\\\text{is decreasing in $(-1,0)$,}
\end{cases}\\[1em]
\text{if $\alpha\ge\dfrac{e}{e-1}$,}&\mu_{\alpha,n}(x)\begin{cases}\text{is
decreasing in $(0,\infty)$,}\\\text{is increasing in $(-1,0)$,}
\end{cases}
\end{cases}
\end{equation*}
as a result, from $\mu_{\alpha,n}(0)=0$, it is easy to obtain that
\begin{equation*}
\mu_{\alpha,n}(x)
\begin{cases}
\ge0&\text{for $\alpha\le0$},\\\le0&\text{for $\alpha\ge\frac{e}{e-1}$},
\end{cases}
\end{equation*}
which is equivalent to
\begin{equation*}
[\ln\nu_\alpha(x)]^{(n)}=\frac{\mu_{\alpha,n}(x)}{x^{n+1}}
\begin{cases}
\ge0&\text{for $\alpha\le0$},\\\le0&\text{for $\alpha\ge\frac{e}{e-1}$}.
\end{cases}
\end{equation*}
In conclusion, we have
\begin{equation}
(-1)^n[\ln\nu_\alpha(x)]^{(n)}\begin{cases} \le0&\text{if
$\alpha\le0$,}\\\ge0&\text{if $\alpha\ge\frac{e}{e-1}$.}
\end{cases}
\end{equation}
\end{rem}

\begin{rem}
It has been proved in the proof of Theorem \ref{gammd} that $\tau(s,t)>0$ for
$(s,t)\in(0,\infty)\times(0,\infty)$. Now we give an upper bound of the function
$\tau(s,t)$ on $(0,\infty)\times(0,\infty)$.
\par
Let $s=\mu t$ for $\mu\in(0,\infty)$. Then we have
\begin{equation}
\label{tmu} \tau(\mu
t,t)=\frac1\mu\biggl[1-\frac{(\mu+1)t+1}{1+t}\Bigl(\frac{t}{1+t}\Bigr)^{\mu
t}\biggr].
\end{equation}
Since the function $\frac{(\mu+1)t+1}{1+t}$ is strictly increasing with
$t\in(0,\infty)$ for fixed $\mu\in(0,\infty)$, it follows that
\begin{equation}
\label{jhg} 1<\frac{(\mu+1)t+1}{1+t}<\mu+1.
\end{equation}
Since the function $\bigl(1+\frac1t\bigr)^t$ is strictly increasing, we see
that
\begin{equation}
\Bigl(\frac{t}{1+t}\Bigr)^{\mu t}=\biggl[\frac{1}{(1+1/t)^t}\biggr]^\mu
\end{equation}
is strictly decreasing with $t\in(0,\infty)$ for fixed $\mu\in(0,\infty)$, therefore
\begin{equation}
\label{jklh} \frac1{e^\mu}<\Bigl(\frac{t}{1+t}\Bigr)^{\mu t}<1.
\end{equation}
Combining \eqref{tmu}, \eqref{jhg} and \eqref{jklh} produces
\begin{equation}
\label{jhk} \tau(\mu
t,t)<\frac1\mu\biggl(1-\frac1{e^\mu}\biggr)<\lim_{\mu\to0}\biggl[\frac1\mu\biggl(1-\frac1{e^\mu}\biggr)\biggr]=1.
\end{equation}
Since $\mu\in(0,\infty)$ and $t\in(0,\infty)$ are arbitrary, so we have $\tau(s,t)<1$
for $(s,t)\in(0,\infty)\times(0,\infty)$.
\par
Recently, the upper bound of $\tau(s,t)$ was improved from $1$ to $\frac13$ in
\cite{2-var-function} and further to $\frac3{10}$ in \cite{2-var-ujevic}.
\end{rem}

\begin{rem}
By definition, it is clear that one of the necessary conditions such that
$\frac{(x+1)^\alpha}{[\Gamma(x+1)]^{1/x}}\in\mathcal{L}[(-1,\infty)]$ is $[\ln
\nu_\alpha(x)]'\ge0$ in $(-1,\infty)$, where $\nu_\alpha(x)$ is defined by
\eqref{nudfn}, which is equivalent to
$$
\alpha\le\frac{(x+1)[{x\psi(x+1)}-{\ln\Gamma(x+1)}]}{x^2}.
$$
Combining this with Theorem~\ref{thmadd} yields
\begin{equation}
\tau_0\ge\frac{x^2}{(x+1)[x\psi(x+1)-\ln\Gamma(x+1)]}-1
\end{equation}
for $x\in(-1,\infty)$.
\par
Straightforward numerical computation by the software M\textsc{athematica} shows that
the maximum $\tau_2$ of $\tau(2,t)$ in $(0,\infty)$ is
\begin{equation}
\tau\left(2,\frac{2+\sqrt7}3\right)=\frac12\left[\frac{2+\sqrt7}3-\frac{\bigl(2+\sqrt7\,\bigr)^3\left(3+\frac{2+\sqrt7}3\right)}{27\Bigl(1+\frac{2+\sqrt7}3\Bigr)^3}\right]=0.264076\dotsm
\end{equation}
and the maximum $\tau_3$ of $\tau(3,t)$ in $(0,\infty)$ is
\begin{equation}
\tau\left(3,\frac59+\frac{\sqrt[3]{2836-54\sqrt{406}}}{18}+\frac{\sqrt[3]{1418+27\sqrt{406}}}{9\sqrt[3]{4}}\right)=0.271807\dotsm.
\end{equation}
\par
If $\alpha\le\frac1{1+\tau_2}=0.79\dotsm$, then $\mu_{\alpha,2}'(x)\le0$ and
$\mu_{\alpha,2}(x)$ decreases in $(-1,\infty)$. Since $\mu_{\alpha,2}(0)=0$
and $\lim_{x\to\infty}[\ln\nu_\alpha(x)]^{(2)}=0$, it is obtained that
$[\ln\nu_\alpha(x)]^{(2)}<0$. Therefore the function
$\nu_\alpha(x)=\frac{[\Gamma(x+1)]^{1/x}}{(x+1)^\alpha}$ is strictly
increasing and logarithmically concave for $\alpha\le\frac1{1+\tau_2}$ in
$(-1,\infty)$. If $\alpha\le\frac1{1+\tau_3}=0.78\dotsm$, then
$\mu_{\alpha,3}'(x)<0$ and $\mu_{\alpha,3}(x)$ decreases in $(-1,0)$ and
$\mu_{\alpha,3}'(x)>0$ and $\mu_{\alpha,3}(x)$ increases in $(0,\infty)$. Thus
$\mu_{\alpha,3}(x)\ge0$ and then $[\ln\nu_\alpha(x)]^{(3)}>0$ in
$(-1,\infty)$. Hence $[\ln\nu_\alpha(x)]^{(2)}$ is strictly increasing in
$(-1,\infty)$ if $\alpha\le\frac1{1+\tau_3}$.
\end{rem}

\begin{rem}
It is proved in \cite{compmon,compmon2} that
$$
\frac{\ln\Gamma(x+1)}x-\ln
x+1=\ln\frac{[\Gamma(x+1)]^{1/x}}x+1\in\mathcal{C}[(0, \infty)]
$$
and tends to $\infty$ as $x\to0$ and to $0$ as $x\to\infty$. A similar result was found in \cite{vogt}: The function
$$
1+\frac{\ln\Gamma(x+1)}x-\ln(x+1)=\ln\frac{[\Gamma(x+1)]^{1/x}}{x+1}+1
$$
belongs to $\mathcal{C}[(-1, \infty)]$ and tends to $1$ as $x\to-1$ and to $0$
as $x\to\infty$. These are special cases of our main results, for examples,
Theorem \ref{gammd} and Theorem \ref{gamme}.
\end{rem}

In what follows, as applications of our main results, we would like to deduce some consequences of the theorems stated in Section \ref{minussec1}.

\begin{prop}\label{logcoms}
The function
\begin{equation}\label{logcdfn}
\frac{[{\Gamma(x+\alpha+1)}]^{1/(x+\alpha)}}{[{\Gamma(x+1)}]^{1/x}}
\end{equation}
belongs to $\mathcal{L}[(-1,\infty)]$ if and only if $\alpha>0$.
\par
For $\alpha\ge1$ and $\beta>0$, the function
\begin{equation}\label{ratcompmon}
\frac{[\Gamma(x+1)]^{1/x}}{[\Gamma(x+1+\beta)]^{1/(x+\beta)}}\biggl(1+\frac{\beta}{x+1}\biggr)^\alpha
\end{equation}
belongs to $\mathcal{L}[(-1,\infty)]$. For $\beta>0$ and any given real number
$\alpha$ satisfying $\alpha\le\frac1{1+\tau_0}<1$, the reciprocal of the
function defined by \eqref{ratcompmon} for $\beta>0$ belongs to
$\mathcal{L}[(-1,\infty)]$.
\par
For $\alpha\ge1$ and $\beta>0$, the function
\begin{equation}\label{ratcompmon1}
\frac{[\Gamma(x+1)]^{1/x}}{[\Gamma(x+1+\beta)]^{1/(x+\beta)}}\biggl(1+\frac{\beta}{x}\biggr)^\alpha
\end{equation}
belongs to $\mathcal{L}[(0,\infty)]$. For $\alpha\le0$ and $\beta>0$, the reciprocal
of the function defined by \eqref{ratcompmon1} belongs to $\mathcal{L}[(0,\infty)]$.
\end{prop}

\begin{proof}
These follow from combining Theorem \ref{ratthm} with Theorem \ref{gammc}, Theorem
\ref{gammd}, Theorem \ref{gamme}, and Theorem \ref{thmadd}.
\end{proof}

\begin{rem}
In \cite{kershawl,minc}, among other things,
 the following monotonicity results were obtained:
\begin{align*}
&\left[\Gamma(1+k)\right]^{1/k}<\left[\Gamma(2+k)\right]^{1/(k+1)},
\quad k\in\tn;\\
&\left[\Gamma\left(1+\frac1x\right)\right]^x\text{ decreases with } x>0.
\end{align*}
These are extended and generalized in \cite{ingamma, ingamma-rgmia,
martnew, new-martin}, among other things: The function
$[\Gamma(r)]^{1/(r-1)}$ is increasing in $r>0$. Clearly, Theorem
\ref{gammc} generalizes these results and extends them for the range
of the argument.
\par
The first conclusion in Proposition \ref{logcoms} shows that the sequences
\begin{equation}
\frac{\sqrt[k]{k!}}{\sqrt[m+k]{(m+k)!}}
\quad\text{and}\quad\frac{\bigl[\sqrt[k]{k!}\,\bigr]
\bigl[\sqrt[k+m+n]{(k+m+n)!}\,\bigr]}{\bigl[\sqrt[k+m]{(k+m)!}\,\bigr]
\bigl[\sqrt[k+n]{(k+n)!}\,\bigr]}
\end{equation}
are increasing with $k\in\tn$ for given natural numbers $m$ and $n$.
\end{rem}

\begin{rem}
The results in Proposition \ref{logcoms} generalize and extend those of
\cite{rsg,indon}.
\end{rem}

Define
\begin{equation}
\label{alzgam}
Q_{a,b}(x)=\frac{[\Gamma(x+a+1)]^{1/(x+a)}}{[\Gamma(x+b+1)]^{1/(x+b)}}
\end{equation}
for nonnegative real numbers $a$ and $b$. J. S\'andor \cite{sandorart} established
that $Q_{1,0}$ is decreasing on $(1,\infty)$. In \cite{alzbergii} Alzer and Berg
proved that $[Q_{a,b}(x)]^c$ is completely monotonic with $x\in(0,\infty)$ if and
only if $a\ge b$ for $c>0$. The following proposition extends the ranges of variables
$a,b$ and $x$ in \cite{alzbergii} and can be regarded as a generalization of
Proposition~\ref{logcoms} above.

\begin{prop}
Let $a,b\in\tr$ and $c>0$. Then $[Q_{a,b}(x)]^c\in\mathcal{L}[(-(1+b),\infty)]$ if
and only if $a>b$.
\end{prop}

\begin{proof}
From Theorem \ref{gammc}, it is clear that
$$
\frac1{[\Gamma(x+a+1)]^{1/(x+a)}}\in\mathcal{L}[(-(1+a),\infty)]
$$
for $a\in\tr$. From Theorem~\ref{ratthm} it follows that the function $Q_{a,b}(x)$ is logarithmically completely monotonic in $(-(1+a),\infty)\cap(-(1+b),\infty)=(-(1+b),\infty)$ for $a>b$. So does the
function $[Q_{a,b}(x)]^c$ for $c>0$.
\par
If $[Q_{a,b}(x)]^c$ is logarithmically completely monotonic for $c>0$, then the
derivative $\{[\ln Q_{a,b}(x)]^c\}'=c[g'(x+a)-g'(x+b)]<0$, where $g(x)$ is defined by
\eqref{gxdef} and $g'(x)$ is strictly decreasing in $(-1,\infty)$, since
$g''(x)=\int_0^1t^2\psi''(xt+1)\td t<0$ from \eqref{gnint} and \eqref{difpsi}.
Therefore, there must be $a>b$.
\end{proof}

\begin{prop}
Let $f$ be a logarithmically completely monotonic function and $g$ a completely
monotonic function. Then the function
$$
f\biggl(a+b\int_\alpha^xg(t)\td t\biggr)
$$
is logarithmically completely monotonic on an interval $I$ if it is defined on $I$, where $b$ is positive and $\alpha\in I$.
\par
In particular, if $f$ is logarithmically completely monotonic, then the following
functions are also logarithmically completely monotonic:
\begin{gather}
f(ax^\alpha+b),\quad \text{where $a$ is nonnegative numbers and $0\le\alpha\le1$,}\\
f\big(a+b\ln(1+x)\big),\quad\text{where $b$ is nonnegative},\\
f(1-e^{-x}),\\
f\big(\arctan\sqrt{x}\,\big).
\end{gather}
\par
If $f(x)$ is completely monotonic on an interval $I$, then the function
$[A-f(x)]^{-\mu}$ is logarithmically completely monotonic on $I$, where $A>f(x)$ for
$x\in I$ and $\mu\ge0$.
\end{prop}

\begin{proof}
These are direct consequences of Theorem \ref{compositelogconv}.
\end{proof}

\begin{rem}
The following are also logarithmically completely monotonic functions:
\begin{gather}
\exp(-ax^\alpha),\quad\text{where $a\ge0$ and $0\le\alpha\le1$},\\
[a+b\ln(1+x)]^{-\mu},\quad\text{where $a\ge0$, $b\ge0$ and $\mu>0$},\\
(a-be^{-x})^{-\mu},\quad\text{where $a\ge b>0$ and $\mu\ge0$}.
\end{gather}
\end{rem}

\begin{rem}
Finally, we pose an open problem: Let $\tau_0=\tau(s_0,t_0)$ be the maximum
value of $\tau(s,t)$ defined by \eqref{taudfn} on the set
$\tn\times(0,\infty)$. Then
$\frac{(x+1)^\alpha}{[\Gamma(x+1)]^{1/x}}\in\mathcal{L}[(-1,\infty)]$ if and
only if $\alpha\le\frac1{1+\tau_0}<1$.
\end{rem}

\end{document}